\documentclass{article}
\usepackage{graphicx,amsmath,amsfonts,amssymb}
\usepackage[usenames,dvipsnames,svgnames]{xcolor} 
\usepackage[margin=1in]{geometry}

\usepackage{fancyhdr}
\newcommand{\helv}{%
    \fontfamily{phv}\fontsize{9}{11}\selectfont}
\usepackage{hyperref} 
\hypersetup{
  urlcolor=black,
  colorlinks=true,
  citecolor=MidnightBlue,
  urlcolor=Bittersweet,
  pdftitle={The Chen-Teboulle algorithm is the proximal point algorithm},
  pdfauthor={Stephen Becker}
}

\newcommand{\<}{\langle}
\renewcommand{\>}{\rangle}

\DeclareMathOperator*{\argmin}{argmin}
\newcommand{\defeq}{\stackrel{\text{\tiny def}}{=}}
\newcommand{\HH}{\mathcal{H}}                  
\newcommand{\st}{\ensuremath{\;\text{such that}\;}}
\DeclareMathOperator{\dom}{dom}

\newcommand{\x}{\vec{x}} 
\newcommand{\A}{\mathcal{A}}
\newcommand{\df}{\partial f}
\newcommand{\dg}{\partial g}
\DeclareMathOperator{\zer}{zero}

\title{The Chen-Teboulle algorithm is the proximal point algorithm}

\author{Stephen Becker~\thanks{University of Colorado Boulder, USA. Work was also performed 2011--2014 while author was at 
IBM Research, Yorktown Heights, NY, USA and at 
Laboratoire Jacques-Louis Lions, University Paris-6, under a fellowship from the 
 Fondation Sciences Mathématiques de Paris  (FSMP) and by a public grant overseen by the French National Research Agency (ANR) as part of the ``Investissements d'Avenir'' program (reference: ANR-10-LABX-0098)}
}

\date{November 22, 2011; posted  \today}

\begin{document}

\maketitle
\begin{abstract}
We revisit the Chen-Teboulle algorithm using recent insights and show that this allows a better bound on the step-size parameter.
\end{abstract}

\thispagestyle{fancy} 

\section{Background}
Recent works such as \cite{HeYuan10} have proposed a very simple yet powerful technique for analyzing optimization methods. The idea consists simply of working with a different norm in the \emph{product} Hilbert space.  We fix an inner product $\< x, y \> $ on $\HH \times \HH^* $. Instead of defining the norm to be the induced norm, we define the primal norm as follows (and this induces the dual norm)
$$ \|x\|_V = \sqrt{ \<Vx,x\>} = \sqrt{ \<x,x\>_V },\quad \|y\|_{V}^* = \|y\|_{V^{-1}} = \sqrt{ \< y,V^{-1}y \> }
= \sqrt{\<y,y\>_{V^{-1}} }$$
for any Hermitian positive definite $V \in \mathcal{B}(\HH,\HH)$; we write this condition as $V \succ 0$.  For finite dimensional spaces $\HH$, this means that $V$ is a positive definite matrix.

We discuss the canonical proximal point method in a general norm; this generality has been known for a long time, and the novelty will be our specific choice of norm. This allows us to re-derive the Chen-Teboulle algorithm~\cite{chen1994proximal}, which, even though it is not widely used, appears to be the first algorithm in a series of algorithms \cite{Chan08,preconditionedADMM,ChambollePock10,HeYuan10,Condat2011,Vu2011}. Among other features, a benefit of these new algorithms is that they can exploit the situation when a function $f$ can be written as $f(x) = h(Ax)$ for a linear operator $A$. In particular, this is useful when the proximity operator~\cite{Moreau1962} of $h$ is easy to compute but the proximity operator of $h \circ A$ is not easy (the prox of $h \circ A$ follows from that of $h$ only in special conditions on $A$; see~\cite{CombettesPesquet07}).

The benefit of this analysis is that it gives intuition, allows one to construct novel methods, simplifies convergence analysis, gives sharp bounds on step-sizes, and extends to product-space formulations easily.

\subsection{Proximal Point algorithm}
All terminology is standard, and we refer to the textbook~\cite{CombettesBook} for standard definitions.
Let $\A$ be a maximal monotone operator, such as a subdifferential of a proper lower semi-continuous convex function, and assume 
$ \zer(\A) \defeq \{ \x : 0 \in \A\x \}$ is non-empty.
The proximal point algorithm is a method for finding some $\x \in \zer(\A)$. It makes use of the fundamental fact:
$$ 0 \in \A \x \quad \iff \quad \tau \x \in \tau \x + \A\x$$
for any $\tau > 0$. This is equivalent to
$$ \x \in ( I + \tau^{-1}\A)^{-1} \x  \defeq J_{\tau^{-1}\A}(\x)$$
where $J$ is the resolvent operator.  Since $A$ is maximal monotone, the resolvent is single-valued and non-expansive, so in fact we look for a fixed point $\x = J_{\tau^{-1}\A}(\x)$. Furthermore, a major result of convex analysis is that the resolvent is firmly non-expansive, which guarantees that the fixed-point algorithm will weakly converge, cf.~\cite[Example\ 5.17]{CombettesBook}, a consequence of the Krasnosel'ski\u{\i} theorem. To be specific, the algorithm is:
$$ \x_{k+1} = (I+\tau^{-1}\A)^{-1} \x_k.$$
There is no limit on the step-size $\tau$ (which is actually allowed to change every iteration) as long as $\tau >0$.

This can be made more general by using the following fact:
$$ 0 \in \A\x \quad \iff \quad V\x \in V\x + \A\x$$
for some $V \succ 0$.
The algorithm is:
$$ \x_{k+1} = (I+V^{-1}\A)^{-1} \x_k  = (V+\A)^{-1}(V\x_k).$$
All the convergence results of the proximal point algorithm still apply, since if $\A$ is maximal monotone in the induced norm on $\HH$, then $V^{-1}\A$ is maximal monotone in the $\|\cdot\|_V$ norm.

\section{Chen-Teboulle algorithm}
For $f,g \in \Gamma_0(\HH)$, consider
\begin{equation} \label{eq:problem}
    \min_{x,z}\; f(x) + g(z) \st Ax=z, \quad\text{or, equivalently,}\quad 
        \min_{x}\; f(x) + g(Ax) 
\end{equation}
along with its Fenchel-Rockafellar dual (see \cite{RockafellarBook,CombettesBook})
\begin{equation}\label{eq:dual}
\min_{v} \; f^*(A^*v) + g^*(-v) 
\end{equation}
where $A$ is a bounded linear operator and $f^*$ is the Legendre-Fenchel conjugate function of $f$. 
We assume strong duality and existence of saddle-points, e.g., $0\in \text{sri}\left( \dom g - A(\dom f) \right)$. 
The necessary and sufficient conditions for the saddle-points (the primal and dual optimal solutions, $(x,v)$) are \cite[Thm.~19.1]{CombettesBook}:
\begin{equation}
0 \in \partial f(x) + A^*\overbrace{\partial g(\underbrace{Ax}_{z})}^{y=-v}, \quad
0\in \overbrace{A\underbrace{\partial f^*(A^*v)}_{x}}^{z} - \partial g^*(-v)
\end{equation}
Therefore it is sufficient to find (letting $v=-y$)
\begin{equation}\label{eq:main}
\boxed{y\in \partial g( z )}, \;\text{with}\;\boxed{z=Ax},\;\text{and}\;
\boxed{0\in \partial f(x) + A^*y}
\end{equation}
(this is consistent with both equations; recall $\partial f^* = (\partial f)^{-1}$ and similarly for $g$~\cite[Cor.~16.2]{CombettesBook}).
%

After our analysis, it will be clear that this can extend to problems such as $f(E_1x+b_1) + g(E_2x+b_2)$. Tseng considers such a case in his Modified Forward-Backward splitting algorithm\cite{Tseng08}. For now, we stay with $Ax=z$ for simplicity.
The Chen-Teboulle method\cite{chen1994proximal} is designed to fully split the problem and avoid any coupled equations involving both $x$ and $z$.

The algorithm proposed is (simplifying the step size $\lambda$ to be constant):
\begin{align}
    p_{k+1} &= y_k + \lambda (Ax_k - z_k ) \quad\text{``predictor" step} \label{eq:CTp}\\
    x_{k+1} &= \argmin f(x)+ \<p_{k+1},Ax\> + \frac{1}{2\lambda}\|x-x_k\|^2
    = (\df + I/\lambda)^{-1}( x_k/\lambda - A^*p_{k+1} )     \\
    z_{k+1} &= \argmin g(z)- \<p_{k+1},z\> + \frac{1}{2\lambda}\|z-z_k\|^2
    = (\dg + I/\lambda)^{-1}( z_k/\lambda +p_{k+1} )     \\
    y_{k+1} &= y_k + \lambda (Ax_{k+1} - z_{k+1} ) \quad\text{``corrector" step}
\end{align}

Convergence to a primal-dual optimal solution is proved for a step-size 
\begin{equation}\label{eq:CTstep-size}
    \lambda < 1/(2L),\quad\text{where}\quad L=\max( \|A\|, 1 ).
\end{equation}
The convergence proof also allows for error in the resolvent computations, provided they are not too large.

\subsection{Scaled norm view-point}
We can recast Eq.~\eqref{eq:main} as $0 \in \A \vec{x}$ where
$$ \A(\x) = \bigg( \df(x) + A^* y, \dg(z) - y, z - Ax \bigg),\quad \x=(x,z,y).$$
For intuition, we write this in the shape of a ``matrix'' operator
\begin{equation}
    \A = \begin{pmatrix} \df & 0 & A^* \\
        0 & \dg & -I \\
        -A & I & 0 \end{pmatrix}
\end{equation}
where ``matrix-multiplication'' is defined $\A \cdot (x,y,z) = \A( \x ) $.

To apply the proximal point algorithm, we must compute $( \tau I + \A)^{-1}$:
\begin{equation}
    (\tau I + \A)^{-1} = \begin{pmatrix} \tau I + \df & 0 & A^* \\
        0 & \tau I + \dg & -I \\
        -A & I & \tau I \end{pmatrix}^{-1}.
\end{equation}
The $x$ and $z$ variables are coupled, so it is not clear how to solve this. Consider now $( V + \A)$, without requiring that $V=\tau I$. Choose a Hermitian and positive-definite $V$ to make the problem block-separable.  There are many potential $V$, so we restrict our attention to $V$ with the same block-structure as $\A$, and we let the diagonal blocks of $V$ be $\tau I$ as in the standard proximal point algorithm.  Now, we choose the upper triangular portion of $V$ to cancel the upper triangular portion of $\A$:
\begin{equation}
    V = \begin{pmatrix} \tau_x I & 0 & -A^* \\
        0 & \tau_z I & I \\
        -A & I & \tau_y I \end{pmatrix},\quad\text{so that}\quad
        (V+\A)^{-1} = \begin{pmatrix} \tau_x I + \df & 0 & 0 \\
        0 & \tau_z I +\dg & 0 \\
        -2A & 2I & \tau_y I \end{pmatrix}^{-1}.
\end{equation}
Thus the computation of $x$ and $z$ is decoupled.  Now, $y$ depends on the current values of $x$ and $z$, but $x$ and $z$ are independent of $y$ so they can be computed first, and the $y$ is updated.
The algorithm is thus:
$$\x_{k+1} = (V + \A)^{-1} V \x_k $$
and in terms of the block coordinates, this is:
\begin{align} \label{eq:newCT}
    x_{k+1} &= (\df + I/\lambda_x)^{-1}( x_k/\lambda_x - A^*y_{k} )    \\
    z_{k+1} &= (\dg + I/\lambda_z)^{-1}( z_k/\lambda_z +y_{k} )     \\
    y_{k+1} &= y_k + \lambda_y (2Ax_{k+1}-Ax_k - 2z_{k+1} +z_k) \label{eq:CTny}
\end{align}
using $\lambda_{\{x,y,z\}}=\tau^{-1}_{\{x,y,z\}}$.
Choosing $\tau_x=\tau_y=\tau_z=1/\lambda$, and re-organizing the steps, we recover the Chen-Teboulle algorithm (the $x$ and $z$ variables correspond exactly, and the $p_{k+1}$ in \eqref{eq:CTp} is the same as the $y_{k+1}$ in \eqref{eq:CTny} ).

\begin{figure}
    \centering
    \includegraphics[width=.8\textwidth]{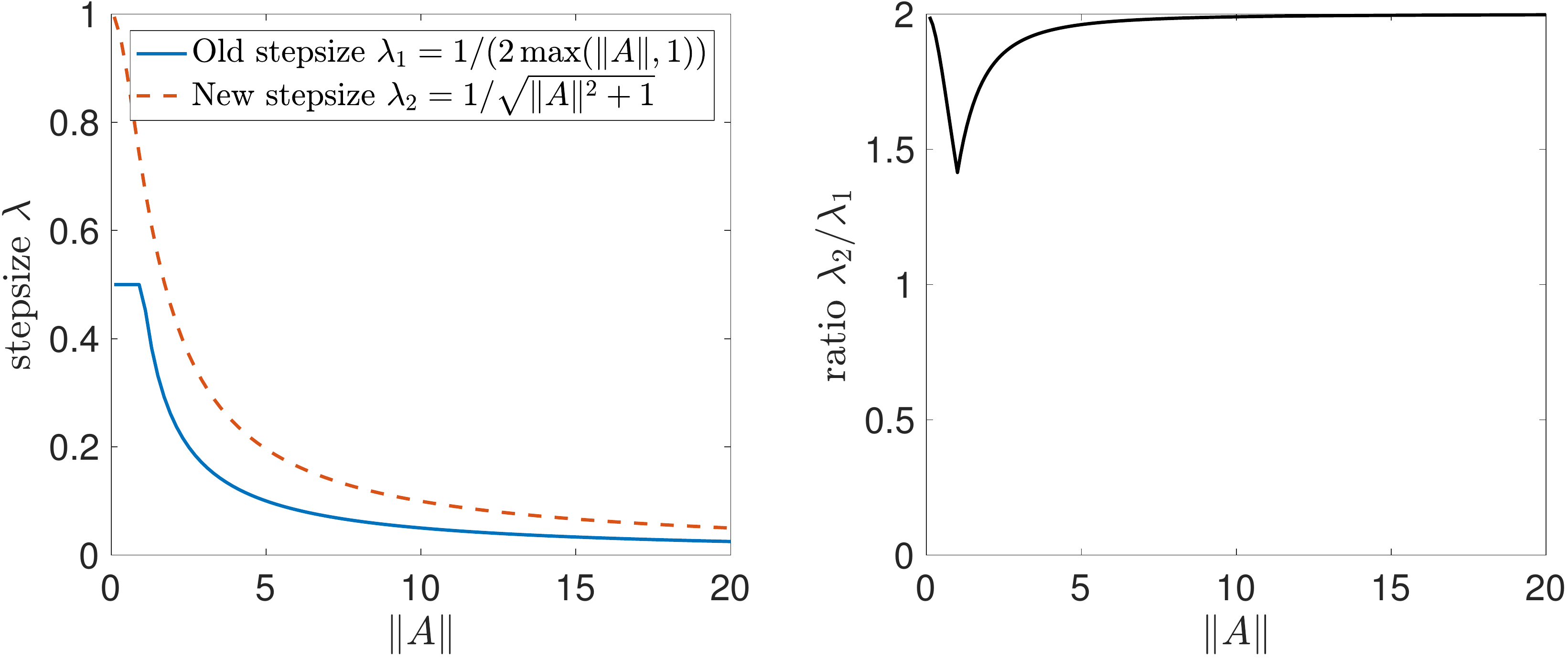}
    \vspace{-1em}
    \caption{The new analysis allows for a larger stepsize}
    \label{fig:1}
\end{figure}

To prove convergence, it only remains to ensure $V \succ 0$. For now, let $V$ be slightly more general:
\begin{equation}
    V = \begin{pmatrix} \tau_x I & 0 & -A^* \\
        0 & \tau_z I & B^* \\
        -A & B & \tau_y I \end{pmatrix}
\end{equation}
where $A$ and $B$ are linear. By applying the Schur complement test twice, we find
$$ V \succ 0 \quad\iff\quad \tau_x > 0,\; \tau_z > 0,\; \tau_y I \succ \tau_x^{-1}AA^* + \tau_z^{-1} BB^*.$$
In the case $B=I, \tau_x=\tau_y=\tau_z=1/\lambda$, then the condition reduces to
\begin{equation}\label{eq:stepsize}
 \lambda \le 1/\sqrt{ \|AA^*\| + 1 } =  1/\sqrt{ \|A\|^2 + 1 } 
 \end{equation}
which is less restrictive than the condition \eqref{eq:CTstep-size} derived in the Chen-Teboulle paper; 
see Fig.~\ref{fig:1}.

An advantage of this approach is that we are free to choose $\tau_x \neq \tau_y \neq \tau_z$. For example, choose $\tau_x = \|AA^*\|$, $\tau_z = \|BB^*\|=1$, $\tau_y < 1/2$.

\subsection*{Acknowledgments}
\vspace{-.5em}
The author is grateful to P.~Combettes for helpful discussions and introducing the Chen-Teboulle algorithm.

\pdfbookmark[1]{References}{refSection}
\bibliographystyle{amsalpha}

\bibliography{thesisBecker}
\end{document}